\documentclass[a4paper]{article}

\usepackage{amsmath,amssymb,amsthm,mathtools}
\usepackage{newtxtext,newtxmath}
\usepackage{bm}
\usepackage{geometry}
\newtheorem{rei}{Example}
\newtheorem{teiri}{Theorem}
\newtheorem{hodai}{Lemma}
\newtheorem{prop}{Proposition}
\newtheorem{remark}{Remark}

\newcommand{\asympL}{\asymp}

\title{Watanabe's Learning Coefficients in Item Response Theory}
\author{Joe Suzuki \and Takumi Ikejiri \and Yuki Kurumadani}

\usepackage{natbib}
\begin{document}
\maketitle

\begin{abstract}
We study the joint-parameter two-parameter logistic (2PL) model in item response theory, in which both item parameters and respondent abilities are treated as unknown fixed parameters. We first characterize its singular structure through the Fisher information matrix. The parameterization has an inherent scale indeterminacy; when all discrimination parameters are nonzero and the ability parameters are not all equal, its tangent direction is the only null direction and the quotient model is locally regular. We then classify the additional rank deficiencies caused by zero discrimination parameters and degenerate ability configurations. This classification separates the strata governed by a nondegenerate quadratic approximation from those requiring an analysis of higher-order product singularities. Under an asymptotic regime in which the numbers of items and respondents are fixed and the number of independent response matrices tends to infinity, we derive explicit local and global learning coefficients and multiplicities for every realizable true response distribution. The global formulas depend only on the true response probabilities, including cases in which some probabilities equal \(1/2\) or are constant across respondents. These results extend previous calculations restricted to nonzero true discrimination parameters and provide the quantities required for applying sBIC to the joint-parameter 2PL model.
\end{abstract}
\noindent\textbf{Keywords:} item response theory; two-parameter logistic model; singular learning theory; learning coefficient; real log canonical threshold; Fisher information matrix; sBIC.

\section{Introduction}
In this paper, we derive the learning coefficient $\lambda$ for a statistical model based on item response theory.

In singular learning theory, established by Sumio Watanabe \citep{Watanabe2011Algebraic}, the learning coefficient is a fundamental quantity that characterizes the asymptotic behavior of the Bayesian generalization error and marginal likelihood \citep{Suzuki2023WAICWBIC}. In particular, it provides the theoretical foundation for the widely applicable information criterion (WAIC) and the widely applicable Bayesian information criterion (WBIC). Conventional information criteria, such as the Akaike information criterion (AIC) and the Bayesian information criterion (BIC), are primarily intended for regular statistical models. By contrast, WAIC and WBIC are applicable to singular, nonregular models. A statistical model $\{p(\cdot\mid\theta)\}$ is said to be regular if the following three conditions hold.
\begin{enumerate}
\item There exists a unique parameter $\theta_*$ that minimizes the Kullback--Leibler divergence from the true distribution $q$.
\item $\theta_*$ lies in the interior of the parameter space $\Theta$.
\item The Fisher information matrix at $\theta_*$ is nonsingular, or equivalently, invertible.
\end{enumerate}
Although the parameter dimension $d$ appears in AIC and BIC, the learning coefficient of a regular model is $\lambda=d/2$. Thus, for a nonregular model, the learning coefficient may be viewed as half its effective number of parameters. One application of the learning coefficient is the information criterion sBIC \citep{DrtonPlummer2017sBIC}. More generally, the learning coefficient governs the asymptotic behavior of the generalization error and marginal likelihood in singular statistical models, provides the theoretical basis for model selection criteria such as WAIC, WBIC, and sBIC, and quantifies effective model complexity.

Several numerical methods for computing learning coefficients have been proposed \citep{Imai2019RLCT}. Explicit formulas, however, are known for only a small number of statistical models. Representative examples include reduced-rank regression \citep{AoyagiWatanabe2005ReducedRankRegression} and analyses of singularities in neural networks \citep{Aoyagi2024LinearNeuralNetworks,Aoyagi2025SingularLearning}. More recently, \citet{Kurumatani2026Semiregular} derived a general formula for the learning coefficient of an extended class of regular models called semiregular models.

The nonregularity of binary item response models has also been noted in other contexts. For example, \citet{ChenMoustakiZhang2020} considered exploratory item factor analysis with the latent factors integrated out and showed that the usual Wilks theory fails when comparing one- and two-factor models because the discrimination parameters associated with the additional factor vanish. Whereas their study concerns comparisons of latent dimensions, the present paper studies the joint-parameter formulation of the one-dimensional 2PL model itself, in which the item parameters and the ability parameter of each respondent are estimated jointly.

For learning coefficients in item response theory, \citet{Ikejiri2024} derived a formula under the assumption that all true discrimination parameters are nonzero. We remove this restriction and derive general formulas that include items with zero discrimination and configurations in which all ability values coincide.\footnote{Throughout this paper, the parameter space for \(a,b,\theta\) is the entire real space, with no sign restrictions imposed.}

Even when all discrimination parameters are nonzero and the ability values are nondegenerate, the joint-parameter 2PL model is not regular in the usual sense because of scale nonidentifiability. At points where some discrimination parameters vanish or all ability values coincide, additional singularities arise beyond the scale direction. Before deriving the learning coefficients, we therefore separate these sources of singularity and determine the null space and rank of the Fisher information matrix.

We consider the two-parameter logistic model of item response theory. For $a_j,b_j,\theta_i\in\mathbb R$, the distribution of $u_{i,j}=\pm1$ is defined by
\[
p(u_{i,j}\mid a_j,b_j,\theta_i)
=
\frac{1}{1+\exp\{-u_{i,j}a_j(\theta_i-b_j)\}}.
\]
Let
\[
u=(u_{ij})_{1\leq i\leq N,\,1\leq j\leq M}.
\]
We assume that the responses are conditionally independent given $(a,b,\theta)$ and define the joint distribution of the response matrix by
\[
p(u\mid a,b,\theta)
=
\prod_{i=1}^{N}\prod_{j=1}^{M}
p(u_{ij}\mid a_j,b_j,\theta_i).
\]
We also let $q$ denote the true joint distribution of the response matrix and $q_{ij}$ its $(i,j)$ marginal distribution, and assume the corresponding product representation
\[
q(u)=\prod_{i=1}^{N}\prod_{j=1}^{M}q_{ij}(u_{ij}).
\]
Here, $M$ is the number of items, $N$ is the number of respondents, and $a=(a_1,\ldots,a_M)$, $b=(b_1,\ldots,b_M)$, and $\theta=(\theta_1,\ldots,\theta_N)$ are the item discrimination parameters, item difficulty parameters, and respondent ability parameters, respectively. We treat $a,b$, and $\theta$ jointly as unknown fixed parameters and do not integrate the ability parameters $\theta$ over an ability distribution. We refer to this formulation as the joint-parameter 2PL model.

Replacing $\theta_i$ and $b_j$ by $\theta_i-\theta_1$ and $b_j-\theta_1$, respectively, leaves the distribution unchanged. We therefore assume that
$$\theta_1=\theta_1^*=0.$$
Accordingly, the parameter space is
\[
\Theta:=\mathbb R^{M}\times\mathbb R^{M}\times\mathbb R^{N-1}.
\]


The main contribution of this paper is to clarify the relationship between the singular structure and the learning coefficients of the joint-parameter 2PL model. First, we analyze the null space of the Fisher information matrix and show that, when all discrimination parameters are nonzero and the ability values are nondegenerate, the only nonidentifiable direction is the one-dimensional scale direction and the quotient model under this action is locally regular. We then classify the additional null directions caused by zero discrimination parameters and degenerate ability values and determine the rank of the Fisher information matrix for every parameter configuration. This classification distinguishes the strata on which the learning coefficient is determined by the quadratic part of the Fisher information matrix from those requiring an analysis of higher-order product singularities. Finally, using resolution of singularities, we derive general formulas for the local learning coefficient and its multiplicity, as well as global formulas expressed solely in terms of the true response probabilities.

Section 2 separates scale nonidentifiability from the additional singularities and classifies the singular strata according to the null space and rank of the Fisher information matrix. Section 3 first develops the general theory of local learning coefficients and multiplicities and then analyzes the basic example in which all true response probabilities equal \(1/2\). Section 4 derives local formulas for arbitrary realizable parameter values and then gives global formulas in terms of the true response probabilities. Section 5 concludes the paper.

\section{Singular Structure of the Joint-Parameter 2PL Model}\label{sec:singular_structure}
In this section, we analyze the singular structure of the joint-parameter 2PL model through its Fisher information matrix. Since the Fisher information matrix represents the quadratic part of the Kullback--Leibler divergence, examining its kernel and rank allows us to distinguish nonidentifiable directions from singular directions that require analysis of higher-order terms. We first isolate the scale indeterminacy that is always present in the model and show that, on the nondegenerate stratum, the quotient model obtained by removing this direction is locally regular. We then classify the additional singularities arising from zero discrimination parameters and degenerate ability parameters.

\subsection{Scale Indeterminacy and Local Regularity of the Quotient Model}
\label{subsec:scale_quotient}

The present \(2M+N-1\)-dimensional parameterization has a scale indeterminacy. At any point at which all discrimination parameters are nonzero and the ability parameters are nondegenerate, the scale direction is the only nonidentifiable direction, and the quotient model obtained by factoring out its action is locally regular.

\begin{prop}\label{prop:quotient_regular}
Let \(P=(a^*,b^*,\theta^*)\in\Theta\), and suppose that \(a_j^*\neq0\) for every \(j\) and that \(\theta_1^*,\ldots,\theta_N^*\) are not all equal. Let \(I(P)\) denote the Fisher information matrix at \(P\), and define
\[
s(P)
=(a_1^*,\ldots,a_M^*,-b_1^*,\ldots,-b_M^*,-\theta_2^*,\ldots,-\theta_N^*)^{\mathsf T}
\]
Then
\[
\ker I(P)=\operatorname{span}\{s(P)\},\qquad
\operatorname{rank}I(P)=2M+N-2
\]
Therefore, the quotient model obtained by factoring out the action of scale transformations is locally regular.
\end{prop}

\begin{proof}
Impose the identifiability constraint \(\theta_1=0\). Let \(J(P)\) denote the Jacobian matrix at \(P\) of the map
\[
\Phi:(a,b,\theta)\longmapsto
\bigl(a_j(\theta_i-b_j)\bigr)_{1\leq i\leq N,\ 1\leq j\leq M}
\]
Writing \(\eta_{ij}=a_j(\theta_i-b_j)\), we have
\[
\left.\frac{\partial\eta_{ij}}{\partial a_k}\right|_P
=\delta_{jk}(\theta_i^*-b_j^*),\qquad
\left.\frac{\partial\eta_{ij}}{\partial b_k}\right|_P
=-\delta_{jk}a_j^*,\qquad
\left.\frac{\partial\eta_{ij}}{\partial\theta_\ell}\right|_P
=\delta_{i\ell}a_j^*
\]
where \(2\leq\ell\leq N\).

For \(u\in\{-1,1\}\),
\[
p(u\mid\eta)=\frac{1}{1+\exp(-u\eta)},\qquad
\frac{\partial}{\partial\eta}\log p(u\mid\eta)
=up(-u\mid\eta)
\]
and
\[
\mathbb E_\eta\left[
\frac{\partial}{\partial\eta}\log p(U\mid\eta)
\right]=0,\qquad
\mathbb E_\eta\left[
\left\{\frac{\partial}{\partial\eta}\log p(U\mid\eta)\right\}^2
\right]
=p(1\mid\eta)p(-1\mid\eta)
\]
hold. Therefore, by the independence of the responses,
\[
I(P)=J(P)^{\mathsf T}W(P)J(P)
\]
where \(\eta_{ij}^*=a_j^*(\theta_i^*-b_j^*)\) and
\[
W(P)=\operatorname{diag}_{i,j}
\left\{p(1\mid\eta_{ij}^*)p(-1\mid\eta_{ij}^*)\right\}
\]
Since \(W(P)\) is positive definite and, for every vector \(x\),
\[
x^{\mathsf T}I(P)x
=\{J(P)x\}^{\mathsf T}W(P)\{J(P)x\}
\]
we have
\[
\ker I(P)=\ker J(P).
\]

For \(c>0\), consider the scale transformation
\[
T_c(a,b,\theta)
=\left(ca,\frac{b}{c},\frac{\theta}{c}\right)
\]
Then
\[
\Phi\{T_c(a,b,\theta)\}=\Phi(a,b,\theta).
\]
Moreover,
\[
\left.\frac{d}{dc}T_c(P)\right|_{c=1}=s(P)
\]
and hence
\[
0
=\left.\frac{d}{dc}\Phi\{T_c(P)\}\right|_{c=1}
=J(P)s(P).
\]
Because \(s(P)\neq0\) by assumption,
\[
\operatorname{span}\{s(P)\}\subseteq\ker I(P)
\]

Conversely, let \((\dot a,\dot b,\dot\theta)\in\ker I(P)\). Since \(\ker I(P)=\ker J(P)\), for every \(i,j\),
\begin{equation}\label{eq:tangent_zero}
(\theta_i^*-b_j^*)\dot a_j
+a_j^*(\dot\theta_i-\dot b_j)=0,
\qquad \dot\theta_1=0
\end{equation}
holds. Substituting \(i=1\) into equation \eqref{eq:tangent_zero} gives
\[
\dot b_j=-\frac{b_j^*}{a_j^*}\dot a_j
\]
Substituting this expression into equation \eqref{eq:tangent_zero} for a general \(i\) yields
\[
\theta_i^*\dot a_j+a_j^*\dot\theta_i=0
\]
Since \(\theta_1^*=0\) and the ability parameters are not all equal, there exists \(r\in\{2,\ldots,N\}\) such that \(\theta_r^*\neq0\). Setting \(i=r\), we find that
\[
\frac{\dot a_j}{a_j^*}
=-\frac{\dot\theta_r}{\theta_r^*}
\]
is constant in \(j\). Denote this value by \(\gamma\). Then
\[
\dot a_j=\gamma a_j^*,\qquad
\dot b_j=-\gamma b_j^*,\qquad
\dot\theta_i=-\gamma\theta_i^*
\]
Therefore,
\[
(\dot a,\dot b,\dot\theta)=\gamma s(P)
\]
and hence
\[
\ker I(P)\subseteq\operatorname{span}\{s(P)\}
\]
Thus,
\[
\ker I(P)=\operatorname{span}\{s(P)\},\qquad
\dim\ker I(P)=1.
\]
Since the parameter space has dimension \(2M+N-1\),
\[
\operatorname{rank}I(P)
=(2M+N-1)-1
=2M+N-2
\]
which proves the result.
\end{proof}

\subsection{Additional Singularities and the Rank of the Fisher Information Matrix}
\label{subsec:singular_strata}

When the nondegeneracy conditions in the preceding subsection fail, kernel directions other than the scale direction arise. These additional singularities fall into two broad categories: those caused by zero discrimination parameters and those caused by all ability parameters being equal. We now classify the rank of the Fisher information matrix and the dimension of its kernel for arbitrary \(P\in\Theta\).

Suppose that at least one discrimination parameter is nonzero and that the ability parameters are not all equal. For an item with \(a_j^*=0\),
\[
\left.\frac{\partial\eta_{ij}}{\partial b_j}\right|_P=-a_j^*=0
\]
so the difficulty direction for that item is added to the kernel.

On the other hand, if all ability parameters are equal, the constraint \(\theta_1^*=0\) implies that
\[
\theta_1^*=\cdots=\theta_N^*=0
\]
For each item with a nonzero discrimination parameter, the linear predictor is then \(-a_j^*b_j^*\), and the item-specific transformation
\[
(a_j,b_j)\longmapsto(c_ja_j,b_j/c_j)
\]
leaves the linear predictor unchanged. Whereas the nondegenerate stratum has only one scale direction common to all items, on the stratum with degenerate ability parameters, independent scale directions arise for the individual items. Moreover, if \(a_j^*=b_j^*=0\) on this stratum, all first derivatives with respect to that item vanish, and, up to a constant factor, the leading part of the Kullback--Leibler divergence contains a product-type term of the form
\[
a_j^2\sum_{i=1}^N(\theta_i-b_j)^2
\]
This higher-order degeneracy cannot be described solely by the rank of the Fisher information matrix.

For arbitrary \(P=(a^*,b^*,\theta^*)\in\Theta\), define
\[
L(P):=\#\{j:a_j^*\neq0\},\qquad
S(P):=\#\{j:a_j^*=0,\ b_j^*=\theta_1^*=\cdots=\theta_N^*\}
\]
When all ability parameters are equal, \(\theta_i^*=0\), and hence \(S(P)\) equals the number of items for which \(a_j^*=b_j^*=0\).

\begin{teiri}\label{thm:fisher_rank}
The rank and the dimension of the kernel of the Fisher information matrix for the joint-parameter 2PL model are given by
\[
\operatorname{rank}I(P)=
\begin{cases}
M+L(P)+N-2,
& L(P)\geq1,\quad
\theta_1^*,\ldots,\theta_N^*\text{ are not all equal},\\[2mm]
M+N-1-S(P),
& L(P)\geq1,\quad
\theta_1^*=\cdots=\theta_N^*,\\[2mm]
M-S(P),
& L(P)=0,
\end{cases}
\]
and, writing the dimension of the parameter space as \(d=2M+N-1\),
\[
\dim\ker I(P)=
\begin{cases}
M-L(P)+1,
& L(P)\geq1,\quad
\theta_1^*,\ldots,\theta_N^*\text{ are not all equal},\\[2mm]
M+S(P),
& L(P)\geq1,\quad
\theta_1^*=\cdots=\theta_N^*,\\[2mm]
M+N-1+S(P),
& L(P)=0.
\end{cases}
\]
\end{teiri}

\begin{proof}
We use equation \eqref{eq:tangent_zero}. First, suppose that \(L(P)\geq1\) and that the ability parameters are not all equal. For each item with \(a_j^*\neq0\), the same argument as in the preceding subsection shows that, for some common \(\gamma\),
\[
\dot a_j=\gamma a_j^*,\qquad
\dot b_j=-\gamma b_j^*,\qquad
\dot\theta_i=-\gamma\theta_i^*
\]
For an item with \(a_j^*=0\), on the other hand, the ability parameters are not all equal, and hence
\[
(\theta_1^*-b_j^*,\ldots,\theta_N^*-b_j^*)\neq0
\]
It follows from equation \eqref{eq:tangent_zero} that \(\dot a_j=0\), whereas \(\dot b_j\) is unrestricted. Thus, the kernel consists of one common scale direction and \(M-L(P)\) difficulty directions, and
\[
\dim\ker I(P)=M-L(P)+1
\]

Next, suppose that \(L(P)\geq1\) and that all ability parameters are equal. Then \(\theta_i^*=0\). Since at least one \(a_j^*\) is nonzero, equation \eqref{eq:tangent_zero} implies that
\[
\dot\theta_i=0,\qquad i=2,\ldots,N
\]
For an item with \(a_j^*\neq0\), \(\dot a_j\) may be chosen freely, after which
\[
\dot b_j=-\frac{b_j^*}{a_j^*}\dot a_j
\]
is determined, yielding one kernel direction. For an item with \(a_j^*=0,b_j^*\neq0\), we have \(\dot a_j=0\), whereas \(\dot b_j\) is unrestricted; for an item with \(a_j^*=b_j^*=0\), both \(\dot a_j\) and \(\dot b_j\) are unrestricted. Therefore,
\[
\dim\ker I(P)
=L(P)+\{M-L(P)-S(P)\}+2S(P)
=M+S(P)
\]

Finally, suppose that \(L(P)=0\). Then equation \eqref{eq:tangent_zero} becomes
\[
(\theta_i^*-b_j^*)\dot a_j=0
\]
All \(\dot b_j\) and \(\dot\theta_2,\ldots,\dot\theta_N\) are unrestricted, while \(\dot a_j\) is unrestricted only for items satisfying \(b_j^*=\theta_i^*\) for every \(i\), that is, for the items counted by \(S(P)\). Therefore,
\[
\dim\ker I(P)
=M+(N-1)+S(P)
=M+N-1+S(P)
\]
Subtracting the dimension of the kernel from \(d=2M+N-1\) in each case yields the rank formulas.
\end{proof}

This rank classification separates, in the joint-parameter 2PL model, the singularities captured by the quadratic approximation of the Kullback--Leibler divergence from the higher-order singularities that it does not capture. On the stratum where all discrimination parameters are nonzero and the ability parameters are nondegenerate, the kernel consists only of a single scale direction, and the quotient model obtained by factoring out its action is locally regular. When a discrimination parameter is zero, a difficulty direction is added to the kernel; when all ability parameters are equal, the common scale direction splits into item-specific scale directions. In particular, at a point where \(S(P)>0\), the first derivatives in both the discrimination and difficulty directions vanish simultaneously, giving rise to a product-type singularity that does not appear in the Fisher information matrix. Therefore, determining the asymptotic complexity in this case requires analysis of the higher-order structure of the Kullback--Leibler divergence.

\section{General Theory of Learning Coefficients and a Basic Example}\label{sec:learning_preliminaries}
This section presents the definitions of the local learning coefficient and multiplicity, together with the basic lemmas needed for the general formula in the next section. We then analyze the case in which all true response probabilities are \(1/2\), as a basic example in which a product-type singularity appears most clearly.

\subsection{Local Learning Coefficient and Multiplicity}

For $(a,b,\theta)\in\Theta$, define the Kullback--Leibler divergence by
\[
K(a,b,\theta)
:=
\sum_{i=1}^{N}
\sum_{j=1}^{M}
\sum_{u\in\{-1,1\}}
q_{ij}(u)
\log
\frac{q_{ij}(u)}
{p(u\mid a_j,b_j,\theta_i)}.
\]
We assume that the true distribution \(q\) is realizable. That is, suppose that there exists
\[
(\bar a,\bar b,\bar\theta)\in\Theta
\]
such that, for all \(i,j\) and \(u\in\{-1,1\}\),
\[
q_{ij}(u)=p(u\mid \bar a_j,\bar b_j,\bar\theta_i).
\]
We then define
\[
\Theta^*:=\{(a,b,\theta)\in\Theta\mid K(a,b,\theta)=0\}.
\]

In the asymptotic setting considered in this paper, \(M\) and \(N\) are fixed. One observation consists of an entire response matrix, and the sample consists of \(n\) independent response matrices generated from the same true joint distribution \(q(u)\), with \(n\to\infty\). Thus, the parameter dimension \(2M+N-1\) remains fixed in this asymptotic limit. This setting differs from the asymptotic framework commonly used in IRT, in which the number of respondents increases as the sample size. Consequently, the results below apply to the fixed-\((M,N)\) setting described here and do not directly provide formulas for the learning coefficient under an asymptotic framework in which the number of respondents increases.

We assume that the prior density \(\varphi(a,b,\theta)\) is bounded and has compact support in the parameter space \(\Theta\). We further assume that the realizing parameters constructed in the proof of Theorem~\ref{thm:global_2pl_logodds}, at which the global learning coefficient and multiplicity stated there are attained, lie in the interior of \(\operatorname{supp}\varphi\), and that \(\varphi\) is positive and analytic in their neighborhoods. These assumptions exclude contributions from infinity and ensure that the global learning coefficient and multiplicity in Theorem~\ref{thm:global_2pl_logodds} are determined by the local singularities considered below.

Fix any \(P\in\Theta^*\). We define the local learning coefficient \(\lambda(P)\) and local multiplicity \(m(P)\) at \(P\) as follows. By Hironaka's resolution of singularities theorem~\cite{Hironaka1964Resolution}, for a sufficiently small neighborhood \(U\) of \(P\), there exist an analytic manifold \(\widetilde U\), covered if necessary by finitely many local charts, and a proper analytic map
\[
g:\widetilde U\longrightarrow U.
\]
Moreover, in the coordinates \(u=(u_1,\ldots,u_d)\) of each local chart, we can write
\begin{equation}\label{eq5}
K(g(u))=
a(u)u_1^{2k_1}\cdots u_d^{2k_d},
\qquad
|g'(u)|=
b(u)|u_1|^{h_1}\cdots |u_d|^{h_d}    
\end{equation}
where \(a(u)\) and \(b(u)\) are analytic functions that do not vanish at the origin, and \(k_1,\ldots,k_d,h_1,\ldots,h_d\) are nonnegative integers. Define the candidate value in this chart by
\[
\lambda_{\mathrm{chart}}=
\min_{1\le i\le d,\ k_i>0}
\frac{h_i+1}{2k_i},
\]
and let \(m_{\mathrm{chart}}\) be the number of indices attaining this minimum. The local learning coefficient \(\lambda(P)\) at \(P\) is defined as the minimum of \(\lambda_{\mathrm{chart}}\) over all charts appearing in a neighborhood of \(g^{-1}(P)\). The value \(m(P)\) is defined as the maximum of \(m_{\mathrm{chart}}\) among the charts satisfying \(\lambda_{\mathrm{chart}}=\lambda(P)\). Thus, the local learning coefficient \(\lambda(P)\) and multiplicity \(m(P)\) are defined at every \(P\in\Theta^*\). The global learning coefficient \(\lambda\) and multiplicity \(m\) are then defined by
\[
\lambda := \min_{P\in\Theta^*}\lambda(P),\qquad
m := \max_{P\in\Theta^*,\ \lambda(P)=\lambda}m(P).
\]

For nonnegative functions \(A\) and \(B\), we write \(A\asymp B\) if there exist positive constants \(c,C\) and a neighborhood \(U\) such that
\[
cB(x)\leq A(x)\leq CB(x),\qquad x\in U.
\]
For example, for \(n\) functions \(A_1,\ldots,A_n\),
$$1+\sum_{i=1}^n A_i(x)^2\asymp 1.$$
Equation~\eqref{eq5} can also be written as
$$
K(g(u))
\asymp u_1^{2k_1}\cdots u_d^{2k_d},
\qquad
|g'(u)|
\asymp|u_1|^{h_1}\cdots |u_d|^{h_d}.  
$$
Moreover, if \(A\asymp B\) and the weight function is positive and analytic in a neighborhood of the origin, then \(A\) and \(B\) have the same local learning coefficient and local multiplicity. In what follows, equivalent nonnegative analytic functions are identified in this sense.

We next state the lemmas used in the proofs in the main text. Lemma~\ref{hodai1} is used to prove that the learning coefficient \(\lambda\) equals \(d/2\) for a regular statistical model.

\begin{hodai}[Positive definite quadratic form]\label{hodai1}
Let \(F\) be a nonnegative analytic function of \(d\) variables defined in a neighborhood of the origin, and suppose that \(F(0)=0\). If the quadratic part \(Q\) of \(F\) at the origin is a positive definite quadratic form, then
\[
 F(z)\asympL Q(z)\asympL \sum_{i=1}^{d}z_i^2
\]
in a neighborhood of the origin. Consequently, the local learning coefficient and local multiplicity of \(F\) at the origin are
\[
 \lambda=\frac{d}{2},\qquad m=1.
\]
\end{hodai}

\begin{proof}
Bounds in terms of the smallest and largest eigenvalues of the symmetric matrix representing \(Q\) yield
\[
 Q(z)\asympL \sum_{i=1}^{d}z_i^2.
\]
Moreover, since \(F\ge0\) and \(F(0)=0\), the linear term at the origin vanishes, and hence
\[
 F(z)=Q(z)+O(\lVert z\rVert^3).
\]
Therefore, \(F(z)\asympL Q(z)\). The last assertion follows by blowing up the origin of \(\sum_{i=1}^{d}z_i^2\).
\end{proof}

The following lemma is also needed to derive the main results of this paper.
\begin{hodai}[\cite{Lin2017}]\label{hodai2}
Let \(K\) be a nonnegative analytic function of the variables \(x_1,\ldots,x_l\). Suppose that there exist
$$B_1,\ldots,B_n \subseteq \{1,\ldots,l\},\qquad B_i\cap B_j=\emptyset\quad (i\neq j),$$
and nonnegative analytic functions \(K_1,\ldots,K_n\), where \(K_i\) depends only on the variables indexed by \(B_i\), such that \(2\leq n\leq l\) and
\[
K(x_1,\ldots,x_l)
\asymp
K_1(\{x_i\}_{i\in B_1})+\cdots+K_n(\{x_i\}_{i\in B_n}).
\]
If the learning coefficient and multiplicity of each \(K_i\) are \((\lambda_i,m_i)\), then the learning coefficient \(\lambda\) and multiplicity \(m\) of \(K\) are
$$\lambda=\sum_{i=1}^n\lambda_i,\qquad
m=\sum_{i=1}^n m_i-n+1.$$
\end{hodai}
When applying Lemma~\ref{hodai2} in what follows, empty variable blocks are omitted. That is, terms for which \(K_i=0\) are not included in the additivity formula. For notational convenience, however, if we formally assign
\[
(\lambda_i,m_i)=(0,1)
\]
to an empty block, the same formula remains valid.

\subsection{The Case in Which All True Response Probabilities Are \(1/2\)}
\label{subsec:half_probability}

Before presenting the general formula, we analyze a basic example of a product-type singularity.

\begin{rei}[The case in which all true response probabilities are \(1/2\)]\label{ex:half_probability}
Let \(P=(0,0,0)\in\Theta^*\). This is equivalent to
\[
q_{ij}(1)=q_{ij}(-1)=\frac12
\]
for all \(i=1,\ldots,N\) and \(j=1,\ldots,M\). Then
\[
\lambda(P)=\frac{M}{2},
\qquad
m(P)=
\begin{cases}
M+1,& N=1,\\[1mm]
1,& N\geq 2.
\end{cases}
\]
\end{rei}

\begin{proof}
By the local equivalence between the Kullback--Leibler divergence for Bernoulli distributions and the squared difference between their log-odds parameters, in a neighborhood of \(P=(0,0,0)\),
\[
K(a,b,\theta)\asympL
\sum_{j=1}^{M}a_j^2\sum_{i=1}^{N}(b_j-\theta_i)^2,
\qquad \theta_1=0.
\]

When \(N=1\),
\[
 K\asympL\sum_{j=1}^{M}a_j^2b_j^2.
\]
Each term has learning coefficient and multiplicity \((1/2,2)\), and the additivity formula for disjoint variable blocks gives
\[
 \lambda(P)=\frac{M}{2},\qquad m(P)=M+1.
\]

Suppose now that \(N\ge2\). For each \(j\),
\[
 b_j^2+\sum_{i=2}^{N}(b_j-\theta_i)^2
\]
is a positive definite quadratic form in \((b_j,\theta_2,\ldots,\theta_N)\). By Lemma~\ref{hodai1},
\[
 K\asympL
 \sum_{j=1}^{M}a_j^2b_j^2
 +\left(\sum_{i=2}^{N}\theta_i^2\right)
  \left(\sum_{j=1}^{M}a_j^2\right).
\]
Blow up the origin of \(a=(a_1,\ldots,a_M)\). In the \(r\)th chart, set
\[
 a_k=a_rc_k\qquad(k\ne r).
\]
The Jacobian is \(|a_r|^{M-1}\), and
\[
 K\asympL a_r^2H_r,
 \qquad
 H_r\asympL
 \sum_{i=2}^{N}\theta_i^2+b_r^2+
 \sum_{k\ne r}c_k^2b_k^2.
\]
The learning coefficient of \(\sum_{i=2}^{N}\theta_i^2+b_r^2\) is \(N/2\), and that of each \(c_k^2b_k^2\) is \(1/2\). Therefore,
\[
 \lambda(H_r)=\frac{N+M-1}{2}.
\]
It follows that
\[
 \lambda(P)=\min\left\{\frac{M}{2},\frac{N+M-1}{2}\right\}=\frac{M}{2},\qquad m(P)=1.
\]
\end{proof}

\section{General Formulas for the Learning Coefficient and Multiplicity}\label{sec:general_results}
In this section, we use the basic example in the preceding section and the existing result for the case in which all true discrimination parameters are nonzero to derive the local learning coefficient and local multiplicity at an arbitrary realizing parameter. We then give the global learning coefficient and multiplicity in terms of the true response probabilities.

\subsection{Local Formulas}

\begin{prop}[\cite{Ikejiri2024}]\label{prop:ikejiri}
Let \(P=(a^*,b^*,\theta^*)\in\Theta^*\). If \(a_j^*\neq0\), \(j=1,\ldots,M\), then the local learning coefficient and local multiplicity at \(P\) are given by
\[
\lambda(P)=
\begin{cases}
0,& M=0,\\[1mm]
\dfrac{M+N-1}{2},
& M\ge1,\ \theta^*_1=\cdots=\theta^*_N,\\[2mm]
\dfrac{2M+N-2}{2},
& M\ge1,\ \theta^*_1,\ldots,\theta^*_N\text{ are not all equal}.
\end{cases},\qquad m(P)=1
\]
\end{prop}
In what follows, we use the notation\footnote{Here, $\#A$ denotes the cardinality of a finite set $A$.}
\[
L(P):=\#\{j: a_j^*\neq 0\},\qquad
S(P):=\#\{j:a_j^*=0,\ b_j^*=\theta_1^*=\cdots=\theta_N^*\},
\]
to extend Proposition \ref{prop:ikejiri}. Proposition \ref{prop:ikejiri} assumes that $L=M$, and hence that $S=0$. We now remove this restriction.
\begin{teiri}\label{thm:local_learning}
For any \(P=(a^*,b^*,\theta^*)\in\Theta^*\), the local learning coefficient and local multiplicity at \(P\) are given by
\[
\lambda(P)=
\begin{cases}
\dfrac{M}{2},
& L(P)=0,\\[2mm]
\dfrac{M+N-1}{2},
& L(P)\ge1,\ \theta_1^*=\cdots=\theta_N^*,\\[2mm]
\dfrac{L(P)+M+N-2}{2},
& L(P)\ge 1,\ \theta_1^*,\ldots,\theta_N^*
\text{ are not all equal},
\end{cases}
\]
\[
m(P)=
\begin{cases}
S(P)+1,& (L(P)=0,\ N=1)\ \text{or}\ 
(L(P)\ge1,\ \theta_1^*=\cdots=\theta_N^*),\\[1mm]
1,& \text{otherwise}.
\end{cases}
\]
\end{teiri}

\begin{proof}
Fix any \(P=(a^*,b^*,\theta^*)\in\Theta^*\). To simplify notation, throughout this proof we write $\lambda$, $m$, $L$, and $S$ instead of $\lambda(P)$, $m(P)$, $L(P)$, and $S(P)$, respectively.

Partition the set $\{1,\ldots,M\}$ as
\[
\begin{aligned}
A_1&:=\{j:a_j^*\neq0\},\\
A_2&:=\{j:a_j^*=0,\ 
(b_j^*-\theta_1^*,\ldots,b_j^*-\theta_N^*)\neq(0,\ldots,0)\},\\
A_3&:=\{j:a_j^*=0,\ 
b_j^*=\theta_1^*=\cdots=\theta_N^*\},
\end{aligned}
\]
where
\[
\# A_1=L,\qquad \# A_2=M-L-S,\qquad \#A_3=S,
\]
and write
$$K(a,b,\theta)\asymp K_1+K_2+K_3,$$
where
$$
K_1:=\sum_{j\in A_1}\sum_{i=1}^N\{a_j(b_j-\theta_i)-a_j^*(b_j^*-\theta_i^*)\}^2,
\qquad 
K_2:=\sum_{j\in A_2}a_j^2\sum_{i=1}^N (b_j-\theta_i)^2,
\qquad
K_3:=
\sum_{j\in A_3}a_j^2\sum_{i=1}^N (b_j-\theta_i)^2 .
$$
Here, if $A_l$ is empty for $l=1,2,3$, we set $K_l=0$. We use the following lemma.
\begin{hodai}\label{hodai3}
\begin{equation}\label{eq201}
K_2\asymp \sum_{j\in A_2}a_j^2  
\end{equation}
If $L\geq 1$, then
\begin{equation}\label{eq202}
K_1+K_3\asymp K_1+ \sum_{j\in A_3}a_j^2b_j^2.  
\end{equation}
\end{hodai}
See Appendix B for the proof.

In view of \eqref{eq201} and \eqref{eq202}, define
$$
K_2'(\{a_j\}_{j\in A_2})=\sum_{j\in A_2}a_j^2,\qquad
K_3'(\{a_j, b_j\}_{j\in A_3})=\sum_{j\in A_3}a_j^2b_j^2.
$$
If $L=0$, then
$$K(a,b,\theta)\asymp K_2'(\{a_j\}_{j\in A_2})
+K_3(\{a_j,b_j\}_{j\in A_3}, \{\theta_i\}_{i=2}^N),
$$
whereas if $L\geq 1$, then
$$K(a,b,\theta)\asymp K_1(\{a_j,b_j\}_{j\in A_1},\{\theta_i\}_{i=2}^N )
+K_2'(\{a_j\}_{j\in A_2})
+K_3'(\{a_j,b_j\}_{j\in A_3}).
$$
Note that the variables $\{b_j\}_{j\in A_2}$ appearing in $K_2$ contribute only a finite constant factor to the local integral and hence, as indicated by $K_2'(\{a_j\}_{j\in A_2})$, do not affect the learning coefficient or multiplicity. Moreover, in both $(K_2',K_3)$ and $(K_1,K_2',K_3')$, the functions involve disjoint sets of variables. Lemma \ref{hodai2} is therefore applicable in both cases.

For $K_1$, we have $a_j^*\neq0$ for every item in $A_1$, and after relabeling the indices, $K_1$ is the same local function as that in Proposition \ref{prop:ikejiri}. Its learning coefficient and multiplicity are therefore obtained by replacing $M$ with $L$ in Proposition \ref{prop:ikejiri}; namely,
\begin{equation}\label{eq11}
\lambda_1=
\begin{cases}
0,& L=0,\\[1mm]
\dfrac{L+N-1}{2},
& L\ge1,\ \theta^*_1=\cdots=\theta^*_N,\\[2mm]
\dfrac{2L+N-2}{2},
& L\ge1,\ \theta^*_1,\ldots,\theta^*_N\text{ are not all equal}.
\end{cases},\qquad m_1=1.
\end{equation}

For $K_2'$, Lemma \ref{hodai1} gives
$$
\lambda_{2}'=\frac{M-L-S}{2},
\qquad
m_{2}'=1.  
$$

If $L=0$, then
\[
K_3
=
\sum_{j\in A_3}a_j^2\sum_{i=1}^N(b_j-\theta_i)^2
\]
is precisely the function in Example \ref{ex:half_probability} with the number of items set to $S=|A_3|$. Hence,
\[
\lambda_3=\frac{S}{2},
\qquad
m_3=
\begin{cases}
S+1,& N=1,\\[1mm]
1,& N\ge2.
\end{cases}
\]
Consequently,
\[
\lambda
=
\lambda_2'+\lambda_3
=
\frac{M-S}{2}+\frac{S}{2}
=
\frac{M}{2},\qquad
m
=
m_2'+m_3-2+1
=
m_3=
\begin{cases}
S+1,& N=1,\\[1mm]
1,& N\ge2.
\end{cases}
\]

If $L\geq 1$, then
$$
K_3^{(j)}(a_j,b_j)=a_j^2b_j^2\quad (j\in A_3)
$$
is already in normal crossing form, and $k=(1,1)$ and $h=(0,0)$ give
\[
\lambda^{(j)}_3=\frac12,\qquad m_3^{(j)}=2.
\]
Moreover, $\displaystyle K_3'=\sum_{j\in A_3}K^{(j)}_3$ depends only on mutually disjoint pairs of variables $(a_j,b_j)$. Hence, Lemma \ref{hodai2} gives
\[
\lambda_{3}'
=
\sum_{j\in A_3}\lambda_3^{(j)}
=
\frac S2,\qquad
m_{3}'
=
\sum_{j\in A_3} m_3^{(j)}-S+1
=
2S-S+1
=
S+1.
\]
Therefore, if $L\ge1$ and $\theta_1^*=\cdots=\theta_N^*$, then \eqref{eq11} yields
\[
\lambda=\lambda_1+\lambda_2'+\lambda_3'=
\frac{L+N-1}{2}
+\frac{M-L-S}{2}
+\frac S2
=
\frac{M+N-1}{2},
\qquad
m=m_1+m_2'+m_3'-2=S+1.
\]
Finally, suppose that \(L \ge 1\) and that \(\theta_1^*,\ldots,\theta_N^*\) are not all equal. Then \(S=0\). Indeed, if \(S>0\), there would exist some \(j\in A_3\) such that
\[
b_j^*=\theta_1^*=\cdots=\theta_N^*,
\]
which would imply that \(\theta_1^*,\ldots,\theta_N^*\) are all equal, contrary to the assumption. Therefore, by \eqref{eq11} and $\lambda_3'=0$, $m_3'=1$,
\[
\lambda=\lambda_1+\lambda_2'+\lambda_3'=
\frac{2L+N-2}{2}
+\frac{M-L}{2}
=
\frac{L+M+N-2}{2},
\qquad
m=m_1+m_2'+m_3'-2=1.
\]
\end{proof}

\begin{remark}
Combining the rank formula in Theorem \ref{thm:fisher_rank} with the formula for the local learning coefficient in Theorem \ref{thm:local_learning}, we obtain
\[
\lambda(P)=\frac{\operatorname{rank}I(P)+S(P)}{2}
\]
for every \(P\in\Theta^*\). Thus, on the stratum where \(S(P)=0\), the local learning coefficient equals one half of the rank of the Fisher information matrix. On a stratum where \(S(P)>0\), by contrast, product-type singularities that are not reflected in the Fisher information matrix make an additional contribution to the local learning coefficient. The multiplicity cannot be determined from the rank of the Fisher information matrix alone.
\end{remark}

\subsection{Global Formulas}

Theorem \ref{thm:local_learning} gives the local learning coefficient and multiplicity \(\lambda(P),m(P)\) at a fixed realizing parameter \(P=(a^*,b^*,\theta^*)\). The global learning coefficient and multiplicity must therefore be obtained as
\[
\lambda=\min_{P\in\Theta^*}\lambda(P),\qquad
m=\max_{\substack{P\in\Theta^*\\ \lambda(P)=\lambda}}m(P).
\]

Let the true correct-response probabilities and their log odds be
\[
p_{ij}^*:=q(u_{ij}=1),\qquad
\ell_{ij}^*:=\log\frac{p_{ij}^*}{1-p_{ij}^*}=a_j^*(\theta_i^*-b_j^*),\qquad
R:=\#\{j:\ell_{1j}^*\neq0\},\qquad
\ell_j^*:=
\begin{pmatrix}
\ell_{1j}^*\\
\vdots\\
\ell_{Nj}^*
\end{pmatrix}.
\]
For $N\geq 2$, also define
\[
D_j:=
\begin{pmatrix}
\ell_{2j}^*-\ell_{1j}^*\\
\vdots\\
\ell_{Nj}^*-\ell_{1j}^*
\end{pmatrix}
=
a_j^*
\begin{pmatrix}
\theta_2^*-\theta_1^*\\
\vdots\\
\theta_N^*-\theta_1^*
\end{pmatrix},
\qquad
D=
{[D_1,\ldots,D_M]},\qquad
L_D:=\#\{j:D_j\neq0\}.
\]
These quantities are determined by the true distribution and do not depend on the choice of realizing parameter.

  \begin{teiri}\label{thm:global_2pl_logodds}
    The global learning coefficient \(\lambda\) and multiplicity \(m\) are given as follows.
    \begin{itemize}
      \item [(1)]
        If \(N=1\), then
        \[
          \lambda=\frac{M}{2},
          \qquad
          m=M-R+1.
        \]
      \item [(2)] 
        If \(N\geq2\), then
        \[
          \lambda=
          \begin{cases}
            \dfrac{M+N-2+L_D}{2},
              & D\neq0,\\[6pt]
            \dfrac{M+N-1}{2},
              & D=0,\ R\geq1,\\[6pt]
            \dfrac{M}{2},
              & D=0,\ R=0,
          \end{cases},\qquad
          m=
          \begin{cases}
            1,
              & D\neq0,\\[4pt]
            M-R+1,
              & D=0,\ R\geq1,\\[4pt]
            1,
              & D=0,\ R=0.
          \end{cases}
        \]
    \end{itemize}
  \end{teiri}
\begin{proof}
In Theorem \ref{thm:local_learning}, $S$ and $L$ depend on the realizing parameter. We first determine the maximum value of $S$ in the cases where it appears in Theorem \ref{thm:local_learning}, namely,
\begin{equation}\label{eq3}
N=1\quad\text{or}\quad N\geq 2,\quad L\geq1,\quad
\theta_1^*=\cdots=\theta_N^*.
\end{equation}
If $\ell_j^*$ is not the zero vector, then $a_j^*\neq0$, so such an item cannot be counted by $S$. Under \eqref{eq3}, an item can be counted by $S$ only if $\ell_{1j}^*=0$, that is,
\begin{equation}\label{eq1}
\displaystyle a_j^*=0\ \text{or}\ b_j^*=\theta_1^*.    
\end{equation}
For each item $j$ satisfying \eqref{eq1}, there exists a corresponding realizing parameterization satisfying
\begin{equation}\label{eq2}
\displaystyle a_j^*=0\ \text{and}\ b_j^*=\theta_1^*.
\end{equation}
Moreover, $S$ is maximized when every $j$ satisfying \eqref{eq1} also satisfies \eqref{eq2}. Thus, in either case in \eqref{eq3}, the maximum value of $S$ is $M-R$.

If $N\geq 2$, $D=0$, and $R\geq 1$, then the second case in \eqref{eq3} applies. Since $L\geq R$, Theorem \ref{thm:local_learning} gives $\lambda=(M+N-1)/2$ in this case. If $N=1$, it gives $\lambda=M/2$.

If $N\geq 2$ and $D\neq0$, then $D_j=0\Longleftrightarrow a_j^*=0$, and hence $L=L_D$ for every realizing parameter. Since $L=0$ would contradict $D\neq0$, Theorem \ref{thm:local_learning} gives
$$ \lambda = \frac{L+M+N-2}{2}=
\frac{L_D+M+N-2}{2},\qquad m=1.$$

Finally, suppose that $N\geq 2$, $D=0$, and $R=0$. Then every \(\ell_j^*\) is the zero vector, and the minimum value $L=0$ can be attained by choosing
\[
a_1^*=\cdots=a_M^*=0.
\]
Theorem \ref{thm:local_learning} therefore gives
$$ \lambda =
\frac{M}{2},\qquad m=1.$$
\end{proof}

\begin{remark}
When \(N=1\), \(R\) is the number of items whose true response probability differs from \(1/2\). When \(N\geq2\), \(D\neq0\) means that the true response probability varies across respondents for at least one item, and \(L_D\) is the number of items exhibiting such variation. In this case, Theorem \ref{thm:global_2pl_logodds} shows that the global learning coefficient depends on \(L_D\), while the multiplicity is one. When \(D=0\), on the other hand, the true response probability for each item is constant across respondents, and \(R\) is the number of items for which this common true response probability differs from \(1/2\). The remaining two cases in Theorem \ref{thm:global_2pl_logodds} are therefore distinguished according to whether \(R\geq1\) or \(R=0\).
\end{remark}

\section{Conclusion}

In this paper, we analyzed the kernel and rank of the Fisher information matrix for the joint-parameter 2PL model and separated the scale nonidentifiability inherent in the model from the additional singularities caused by zero discrimination parameters and degenerate ability configurations. When all discrimination parameters are nonzero and the ability parameters are not all equal, the only nonidentifiable direction is the single scale direction, and the quotient model obtained by factoring out its action is locally regular. By contrast, zero discrimination parameters or equal ability parameters give rise to additional null directions. In particular, we showed that at points where \(S(P)>0\), product singularities arise that are not visible in the Fisher information matrix.

Furthermore, under an asymptotic regime in which \(M,N\) are fixed and the number of independent response matrices increases, we derived the local learning coefficient and local multiplicity at every realizable parameter and obtained formulas for the global learning coefficient and multiplicity expressed solely in terms of the true response probabilities. In particular,
\[
\lambda(P)=\frac{\operatorname{rank}I(P)+S(P)}{2}
\]
holds. On strata where \(S(P)=0\), the local learning coefficient equals one-half the rank of the Fisher information matrix, whereas on strata where \(S(P)>0\), higher-order product singularities make an additional contribution. These results extend previous results restricted to the case in which all true discrimination parameters are nonzero and provide the theoretical quantities required for applying sBIC to the joint-parameter 2PL model.

This paper focused on the 2PL model. An important direction for future research is to extend the analysis to the 3PL model, which includes guessing parameters. The 3PL model is given by
\[
\Pr(u_{ij}=1\mid a_j,b_j,\theta_i,c_j)
=
c_j+(1-c_j)\frac{1}{1+\exp\{-a_j(\theta_i-b_j)\}},
\qquad 0<c_j<1,
\]
and the additional parameters may give rise to new forms of nonidentifiability and new singular strata. Explicitly determining the learning coefficients and multiplicities of this model is an important topic for future research.

\bibliography{ref}

\section*{Appendix A: Derivation of Proposition \ref{prop:ikejiri}}

\begin{proof}
The case $M=0$ is immediate. Henceforth, suppose that $M\ge1$. Let
\[
 x_j:=a_jb_j-a_j^*b_j^*,\qquad
 y_{ij}:=a_j\theta_i-a_j^*\theta_i^*\quad(i=2,\ldots,N).
\]
Since \(\theta_1=\theta_1^*=0\), the local equivalence between the Kullback--Leibler divergence for Bernoulli distributions and the squared difference between their log-odds parameters gives
\[
K(a,b,\theta)\asymp
\sum_{j=1}^{M}x_j^2+
\sum_{i=2}^{N}\sum_{j=1}^{M}(x_j-y_{ij})^2.
\]
The right-hand side is a positive definite quadratic form in the variables \((x,y)\). Therefore, by Lemma \ref{hodai1} and the preceding local equivalence,
\begin{equation}\label{eq44}
K(a,b,\theta)\asymp \sum_{j=1}^Mx_j^2+\sum_{i=2}^N\sum_{j=1}^My_{ij}^2.
\end{equation}
Moreover, since $a_j\ne0$ in a neighborhood of the true parameter, $(a,b,\theta)\mapsto(a,x,\theta)$ is a regular local change of coordinates.

First, suppose that $\theta_1^*=\cdots=\theta_N^*$. Since $\theta_1^*=0$, all $\theta_i^*$ are zero. In this case, $y_{ij}=a_j\theta_i$, and since $\sum_j a_j^2\asympL1$ in a neighborhood of the true parameter,
\[
 K(a,b,\theta)\asympL
 \sum_{j=1}^{M}x_j^2+
 \sum_{i=2}^{N}\theta_i^2.
\]
Therefore,
\[
 \lambda(P)=\frac{M+N-1}{2},\qquad m(P)=1.
\]

Next, suppose that $\theta_1^*,\ldots,\theta_N^*$ are not all equal. Then $\theta_r^*\ne0$ for some $r\in\{2,\ldots,N\}$.

Starting from the local coordinates $(a,x,\theta)$ used above, we make the following change of variables:
$$ u_j:=a_j\theta_r-a_j^*\theta_r^*\quad(j=1,\ldots,M),
 \qquad
 v_i:=\frac{\theta_i}{\theta_r}-\frac{\theta_i^*}{\theta_r^*}
 \quad(i\ne1,r).    $$
The transformed coordinates are
\[
 (x_1,\ldots,x_M,\,u_1,\ldots,u_M,\,(v_i)_{i\ne1,r},\,\theta_r).
\]
The inverse transformation is given by
\[
 a_j=\frac{u_j+a_j^*\theta_r^*}{\theta_r},
 \qquad
 \theta_i=\theta_r\left(v_i+\frac{\theta_i^*}{\theta_r^*}\right)
 \quad(i\ne1,r).
\]
Since $\theta_r^*\ne0$, this is a regular local change of coordinates in a neighborhood of the true parameter.

Separating the second term in \eqref{eq44} into the cases $i=r$ and $i\ne r$ gives
$$
\begin{aligned}
 \sum_{i=2}^{N}\sum_{j=1}^{M}
 (a_j\theta_i-a_j^*\theta_i^*)^2
 =
 \sum_{j=1}^{M}u_j^2
 +\sum_{\substack{i=2\\i\ne r}}^{N}\sum_{j=1}^{M}
 \left\{
 \frac{\theta_i^*}{\theta_r^*}u_j
 +(u_j+a_j^*\theta_r^*)v_i
 \right\}^2.
\end{aligned}
$$
Its quadratic part is
\[
 \sum_{j=1}^{M}u_j^2
 +\sum_{\substack{i=2\\i\ne r}}^{N}\sum_{j=1}^{M}
 \left(
 \frac{\theta_i^*}{\theta_r^*}u_j
 +a_j^*\theta_r^*v_i
 \right)^2.
\]
If this is zero, the first sum implies that all $u_j=0$, and $a_j^*\theta_r^*\ne0$ then implies that all $v_i=0$. Hence, this quadratic part is positive definite. By Lemma \ref{hodai1},
\begin{equation}\label{eq77}
 \sum_{i=2}^{N}\sum_{j=1}^{M}
 (a_j\theta_i-a_j^*\theta_i^*)^2
 \asympL
 \sum_{j=1}^{M}u_j^2
 +\sum_{\substack{i=2\\i\ne r}}^{N}v_i^2.
\end{equation}
It follows from \eqref{eq44} and \eqref{eq77} that
\[
 K(a,b,\theta)\asympL
 \sum_{j=1}^{M}x_j^2
 +\sum_{j=1}^{M}u_j^2
 +\sum_{\substack{i=2\\i\ne r}}^{N}v_i^2.
\]
The right-hand side is a sum of squares in $2M+N-2$ variables, and the remaining transformed coordinate $\theta_r$ does not appear in this expression. Therefore,
\[
 \lambda(P)=\frac{2M+N-2}{2},\qquad m(P)=1.
\]
\end{proof}

\section*{Appendix B: Proof of Lemma \ref{hodai3}}

\begin{proof}
For $j\in A_2$,
\[
 \sum_{i=1}^{N}(b_j^*-\theta_i^*)^2>0.
\]
Hence, in a neighborhood of the true parameter,
\[
 \sum_{i=1}^{N}(b_j-\theta_i)^2\asympL1.
\]
Therefore,
\[
 K_2
 =\sum_{j\in A_2}a_j^2\sum_{i=1}^{N}(b_j-\theta_i)^2
 \asympL\sum_{j\in A_2}a_j^2.
\]

Next, suppose that $L\ge1$. The case $A_3=\varnothing$ is immediate, so suppose that $A_3\ne\varnothing$. Then, by the definition of $A_3$ and the condition $\theta_1^*=0$,
\[
 \theta_1^*=\cdots=\theta_N^*=0.
\]

For $j\in A_1$, let
\[
 x_j:=a_jb_j-a_j^*b_j^*.
\]
Then
\[
 K_1
 =\sum_{j\in A_1}
 \left\{
 x_j^2+\sum_{i=2}^{N}(x_j-a_j\theta_i)^2
 \right\}.
\]
Therefore,
\begin{align}\label{eq110}
 K_1+K_3
 =\sum_{j\in A_1}
 \left\{
 x_j^2+\sum_{i=2}^{N}(x_j-a_j\theta_i)^2
 \right\}
 +\sum_{j\in A_3}
 \left\{
 (a_jb_j)^2+\sum_{i=2}^{N}(a_jb_j-a_j\theta_i)^2
 \right\},
\end{align}
whereas
\begin{align}\label{eq111}
 K_1+\sum_{j\in A_3}a_j^2b_j^2
 =\sum_{j\in A_1}
 \left\{
 x_j^2+\sum_{i=2}^{N}(x_j-a_j\theta_i)^2
 \right\}
 +\sum_{j\in A_3}(a_jb_j)^2.
\end{align}

Since $a_j\ne0\ (j\in A_1)$ in a neighborhood of the true parameter, \eqref{eq110} and \eqref{eq111} are both positive definite quadratic forms in
\begin{equation}
 \{x_j\}_{j\in A_1},\qquad
 \{a_jb_j\}_{j\in A_3},\qquad
 \theta_2,\ldots,\theta_N.
\end{equation}
Therefore, by Lemma \ref{hodai1},
\[
 K_1+K_3
 \asympL
 \sum_{j\in A_1}x_j^2
 +\sum_{i=2}^{N}\theta_i^2
 +\sum_{j\in A_3}(a_jb_j)^2
 \asympL
 K_1+\sum_{j\in A_3}a_j^2b_j^2.
\]
This proves Lemma \ref{hodai3}.
\end{proof}
\end{document}